\documentclass[oneside,a4paper,11-pt,notitlepage]{article}
\usepackage[left=0.8in, right=0.8in, top=1.5in, bottom=1.5in]{geometry}

\usepackage[T1]{fontenc} 
\usepackage[utf8]{inputenc} 
\usepackage[english]{babel}
\usepackage{lipsum} 
\usepackage{lmodern}
\usepackage{amssymb}
\usepackage{amsthm}
\usepackage{bm}
\usepackage{mathtools}
\usepackage{braket}
\usepackage{esint}
\newcommand{\abs}[1]{{\left|#1\right|}}
\newcommand{\norma}[1]{{\left\Vert#1\right\Vert}}

\usepackage{booktabs}
\usepackage{graphicx}
\usepackage{tikz}
\usetikzlibrary{patterns}
\usepackage{multicol}
\usepackage{caption}
\usepackage{enumerate}
\usepackage[skins,theorems]{tcolorbox}
\tcbset{highlight math style={enhanced,
		colframe=black,colback=white,arc=0pt,boxrule=1pt}}
\def\XXint#1#2#3{{\setbox0=\hbox{$#1{#2#3}{\int}$}
    \vcenter{\hbox{$#2#3$}}\kern-.5\wd0}}

\theoremstyle{definition}
\newtheorem{definizione}{Definition}[section]
\theoremstyle{plain}
\newtheorem{teorema}{Theorem}[section]

\newtheorem{prop}[teorema]{Proposition}

\theoremstyle{definition}
\newtheorem{esempio}{Example}[section]
\newtheorem{oss}[esempio]{Remark}
\newtheorem{open}[esempio]{Open problem}

\DeclareMathOperator{\R}{\mathbb{R}}

\DeclareMathOperator{\diam}{\, \textup{diam}}

\makeatletter
\newcommand{\myfootnote}[2]{\begingroup
	\def\@makefnmark{}%
	\addtocounter{footnote}{-1}%
	\footnote{\textbf{#1} #2}
	\endgroup}
\makeatother

\usepackage{hyperref}
\hypersetup{linktoc=none, bookmarksnumbered, colorlinks=true, linkcolor=red}

 \title{Sharp and quantitative estimates for the $p-$Torsion of convex sets}
\author{Amato V., Masiello A. L., Paoli G.*, Sannipoli R. }
\date{}

\newcommand{\Addresses}{{
  \bigskip 
  \footnotesize 
 
   \textit{E-mail address}, G.~Paoli (Corresponding author)*: \texttt{gloria.paoli@fau.de} 
   
 \medskip
 
 \textsc{ Department of Data Science (DDS)
Chair in Dynamics, Control and Numerics (Alexander von Humboldt-Professorship),
Cauerstr. 11,
91058 Erlangen, Germany.}

  \medskip 
 
  \textit{E-mail address}, V.~Amato: \texttt{vincenzo.amato@unina.it} 
 
   \medskip 
 
  \textit{E-mail address}, A.L.~Masiello: \texttt{albalia.masiello@unina.it} 
   \medskip 
 
  \textit{E-mail address}, R.~Sannipoli: \texttt{rossano.sannipoli@unina.it} 
  
     \medskip 
     
    \textsc{Dipartimento di Matematica e Applicazioni ``R. Caccioppoli'', Universit\`a degli studi di Napoli Federico II, Via Cintia, Complesso Universitario Monte S. Angelo, 80126 Napoli, Italy.}\par\nopagebreak 

}}

\begin{document}

\maketitle

\begin{abstract} 
  Let $\Omega\subset\mathbb{R}^n$, $n\geq 2$, be a non-empty, bounded, open and convex set and let $f$ be a positive and non-increasing function depending only on the distance from the boundary of $\Omega$. We consider the $p-$torsional rigidity associated to $\Omega$ for the Poisson problem with Dirichlet boundary conditions, denoted by $T_{f,p}(\Omega)$. Firstly,  we prove a P\'olya type lower bound for $T_{f,p}(\Omega)$ in any dimension; then, we consider the planar case and  we provide two quantitative estimates  in the case $f\equiv 1 $.      \\


\end{abstract}

\section{Introduction}
Let $\Omega\subset\mathbb{R}^n$, $n\geq 2$, be a non-empty, bounded, open and convex set and let $p\in(1,+\infty)$. We consider the Poisson equation for the $p-$Laplace operator, defined as  
\begin{equation*}
	-\Delta_p u:=-{\rm div}\left(|\nabla  u|^{p-2}\nabla  u\right),
	\end{equation*} 
	with Dirichlet boundary condition: 
\begin{equation}
\label{first_problem}
    \begin{cases}
    -\Delta_p u(x) = f(d(x,\partial \Omega)) & \text{ in  } \Omega\\
    u=0  & \text{ on } \partial\Omega,
    \end{cases}
\end{equation}
 where   $f: [0 , R_\Omega ] \to [0, +\infty[$ is  a continuous, non-increasing and not identically zero function,  $ d(\cdot, \partial \Omega):\Omega \to [0,+\infty[$ is the distance function from the boundary defined as
 \begin{equation}
     \label{disti}
     d(x,\partial\Omega):=\inf_{y\in\partial\Omega}\abs{x-y}
 \end{equation}
  and $R_\Omega$ is the inradius of $\Omega$, i.e.
 \begin{equation}
     \label{inradiuss}
     R_\Omega=\sup_{x\in \Omega} d(x,\partial\Omega).
 \end{equation}
  This class of functions, depending only on the distance, are the so called web functions, see as a reference \cite{gazzola_survay}. 
A function $u\in W^{1,p}_0(\Omega)$ is a weak solution to \eqref{first_problem} if and only if
\begin{equation*}
     \int_\Omega |\nabla u(x)|^{p-2} \nabla u(x) \nabla \varphi(x) \, dx=\int_\Omega f(d(x,\partial\Omega))\varphi(x ) \, dx \quad \forall \varphi \in W^{1,p}_0(\Omega). 
\end{equation*} 
The $(f,p)$-torsional rigidity of $\Omega$, that we denote by $T_{f,p}(\Omega)$,  is defined as
\begin{equation}
\label{variational}
    T_{f,p}(\Omega) = \max_{\substack{\varphi\in W^{1,p}_0(\Omega)\\ \varphi\not \equiv0} }  \frac{\displaystyle{\left(\int_\Omega f(d(x,\partial \Omega))\abs{\varphi(x)} \, dx\right)^{\frac{p}{p-1}}}}{\left(\displaystyle{\int_\Omega \abs{\nabla \varphi(x)}^p \, dx}\right)^{\frac{1}{p-1}}}
\end{equation}
and, if $u_p\in W^{1,p}_0(\Omega)$ is the unique solution to \eqref{first_problem}, we have
\begin{equation*}
     T_{f,p}(\Omega)= \int_\Omega f u_p \, dx. 
 \end{equation*}
For the sake of simplicity, when $f\equiv 1$ in $\Omega$, we set $T_p(\Omega):=T_{1,p}(\Omega)$ and, if we are also in the case $p=2 $, we set $T(\Omega):=T_{1,2}(\Omega)$. We recall that the quantities $T(\Omega)$ and $T_p(\Omega)$ are usually called, respectively, torsional rigidity and $p-$torsional rigidity and so, by analogy,  we have chosen the above terminology for $T_{f,p}(\Omega)$. 

   In what follows, we denote by $|\Omega|$ and $P(\Omega)$ respectively the Lebesgue measure and the perimeter of $\Omega$ in the sense of De Giorgi. 
 In \cite{polya1960} the author gives some estimates on  the torsional rigidity $T(\Omega)$.
In particular, he proves that,  among all bounded, open and convex planar sets, the following inequality holds
\begin{equation}
    \label{pol}
    \frac{T(\Omega)P^2(\Omega)}{\abs{\Omega}^3}\ge \frac{1}{3}
\end{equation}
and equality is asymptotically achieved by a sequence of thinning rectangles. Moreover,  Makai in \cite{makai} proves that among all bounded, open and convex planar sets, the following upper bound holds
\begin{equation}\label{makai_intro}
   \frac{T(\Omega)P^2(\Omega)}{\abs{\Omega}^3}\le \frac{2}{3}, 
\end{equation}
which is sharp on a sequence of thinning triangles (for the exact definition of thinning domains see Definition \ref{thin_rect}).
Estimates \eqref{pol} and \eqref{makai_intro} are generalized to the $p-$Laplacian in \cite{fragala_gazzola_lamboley}. More precisely, the authors prove that, in the class of bounded, open and convex planar sets,
\begin{equation}
    \label{fgl}
     \frac{1}{q+1}< \frac{T_p(\Omega) P^q(\Omega)}{\abs{\Omega}^{q+1}}< \frac{2^{q+1}}{(q+2)(q+1)}, \qquad q=\frac{p}{p-1},
\end{equation}
where the lower and the upper bounds hold asymptotically on a sequence of thinning rectangles and on a sequence of thinning isosceles triangles, respectively. In \cite{gavitone_2014} the authors generalize the lower bound \eqref{fgl} in every dimensions, proving that for open, bounded and convex sets $\Omega\subseteq\mathbb{R}^n$, it holds
\begin{equation}\label{gavit}
 \frac{T_p(\Omega) P^q(\Omega)}{\abs{\Omega}^{q+1}}> \frac{1}{q+1},
\end{equation}
and they extend this result also to the anisotropic case.

 We also recall that in \cite{buttazzo2020convex} the authors consider the functional 
$$ H_k(\Omega)=\dfrac{P(\Omega) T^{k}(\Omega)}{|\Omega|^{\alpha_k}}, \qquad \alpha_k=1+k+\frac{2k-1}{n},$$
and prove that, among bounded, open and convex sets in $\mathbb{R}^n$, this functional  is bounded if and only if $k=1/2$. More precisely, they prove the following:
\begin{equation}\label{butta}
    \frac{1}{\sqrt{3}}\le H_{\frac{1}{2}}(\Omega) \le \frac{2^n n^{3n/2}}{\omega_n }\left(\frac{n}{n+2}\right)^{\frac{1}{2}},
\end{equation}
where $\omega_n$ is the Lebesgue measure of the unit ball. We note that, in the planar case, the lower bound in \eqref{butta} coincides with the one given in \eqref{pol}, while the upper bound is strictly larger than the one given in \eqref{makai_intro}. It is conjectured that, in the higher dimensional case, the upper bound is
$$ H_{\frac{1}{2}}(\Omega)\le n \left(\frac{2}{(n+1)(n+2)}\right)^{\frac{1}{2}}. $$
Moreover, we observe that the lower bound in \eqref{butta} is asymptotically achieved by a sequence of thinning cylinders. 
More precisely, denoting by $w_\Omega$ the minimal width and by $\diam(\Omega)$ the diameter of the set (see Section \ref{section_notion} for the exact definitions), we give the following

\begin{definizione}
    Let $\Omega_l$ be a sequence of non-empty, bounded, open and convex sets of $\mathbb{R}^n$. We say that $\Omega_l$ is a sequence of thinning domains if
    \begin{equation}
        \dfrac{w_{\Omega_l}}{\diam(\Omega_l)}\xrightarrow{l \to 0}0.
    \end{equation}
 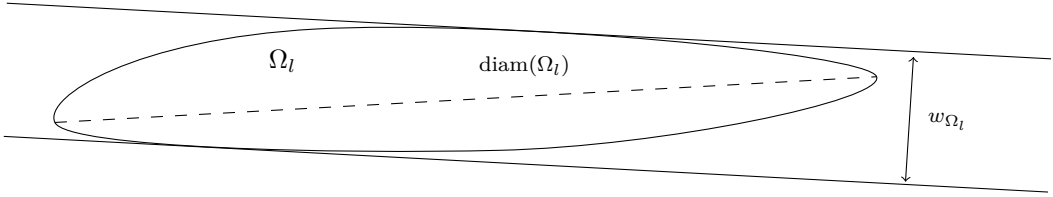
\begin{figure}[h!]
\begin{center}
 \begin{tikzpicture}[x=0.75pt,y=0.75pt,yscale=-1,xscale=1]
\draw   (147.5,317) .. controls (141.5,304) and (187.5,275) .. (277.5,270) .. controls (367.5,265) and (553.5,282) .. (557.5,294) .. controls (561.5,306) and (450.5,329) .. (369.5,331) .. controls (288.5,333) and (153.5,330) .. (147.5,317) -- cycle ;
\draw    (123.5,257) -- (644.5,285) ;
\draw    (122,324) -- (643,352) ;
\draw   [<->]  (575.87,284.05) -- (571.9,346.79) ;
\draw  [dash pattern={on 4.5pt off 4.5pt}]  (147.5,317) -- (557.5,294) ;

\draw (356,282) node [anchor=north west][inner sep=0.75pt]   [align=left] {\footnotesize $\diam(\Omega_l)$};
\draw (582,311) node [anchor=north west][inner sep=0.75pt]   [align=left] {\footnotesize $w_{\Omega_l}$};
\draw (253,279) node [anchor=north west][inner sep=0.75pt]   [align=left] {$\Omega_l$};
\end{tikzpicture}
\end{center}
\caption{Minimal width and diameter of a convex set.} \label{fig:M1}
\end{figure}

    In particular, if $l>0$ and  $C$ is a bounded, open and convex set of $\R^{n-1}$ with unitary $(n-1)$-dimensional measure, then, if $l \to 0$, the sequence 
    \begin{equation}\label{thin_rect}
    \Omega_l = l^{-\frac{1}{n-1}}C \times \left[-\frac{l}{2}, \frac{l}{2}\right]
\end{equation}
    is called a sequence of thinning cylinders. Moreover, in the case $n=2$, the sequence \eqref{thin_rect} is called sequence of thinning rectangles.
    \begin{figure}[h]
       \begin{center}
           \begin{tikzpicture}[x=0.75pt,y=0.75pt,yscale=-1,xscale=1]

\draw  [fill={rgb, 255:red, 230; green, 230; blue, 230 }  ,fill opacity=0.7 ] (237.5,251) .. controls (276.5,220) and (313,212) .. (361.5,236) .. controls (410,260) and (432,275) .. (374.5,297) .. controls (317,319) and (219.5,321) .. (207.5,302) .. controls (195.5,283) and (198.5,282) .. (237.5,251) -- cycle ;
\draw   (237.5,204) .. controls (276.5,173) and (313,165) .. (361.5,189) .. controls (410,213) and (432,228) .. (374.5,250) .. controls (317,272) and (219.5,274) .. (207.5,255) .. controls (195.5,236) and (198.5,235) .. (237.5,204) -- cycle ;
\draw  [dash pattern={on 4.5pt off 4.5pt}]  (304.5,175) -- (304.5,222) ;
\draw    (201.5,240) -- (201.5,287) ;
\draw    (408.5,224) -- (408.5,271) ;
\draw  [<->]  (425.5,226) -- (425.5,273) ;

\draw (280,290) node [anchor=north west][inner sep=0.75pt]   [align=left] {\footnotesize $l^{-\frac{1}{n-1}}C$};
\draw (433,239) node [anchor=north west][inner sep=0.75pt]   [align=left] {\footnotesize $l$};
\draw (200,180) node [anchor=north west][inner sep=0.75pt]   [align=left] { $\Omega_l$};

\end{tikzpicture}
\end{center}
\caption{Thinning cylinders.} \label{fig:M2}
\end{figure}
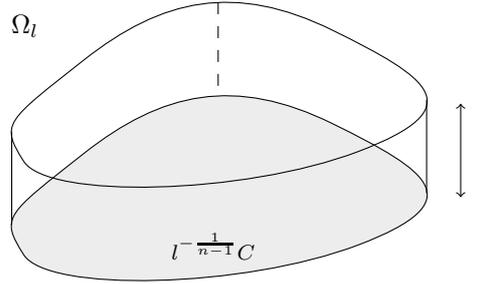
\end{definizione}

The first result that we prove  is a lower bound for the $(f,p)$-torsional rigidity, which generalizes the lower bound in \eqref{fgl}. The hypothesis that $f$ is a web function allows us to use the method of proof contained in \cite{polya1960}.


\begin{teorema}
\label{teorema1}
Let $\Omega$ be a non-empty, bounded, open and convex set of $\R^n$, $n\ge2$, and let $f: [0 , R_\Omega ] \to [0, +\infty[$ be   a continuous and non-increasing function such that  $f\not\equiv 0$. Then, 
\begin{equation}
    \label{polyatype}
    T_{f,p}(\Omega)\ge c_p\frac{\mu^{q+1}_f(\Omega)}{f(0)P^{q}(\Omega)},
\end{equation}
where 
$$ c_p =\frac{p-1}{2p-1}, \qquad q=\frac{p}{p-1},$$
and
$$
    \mu_f(\Omega)= \int_\Omega f(d(x,\partial \Omega)) \, dx.
$$
Moreover, the equality sign is asymptotically achieved by a sequence of thinning cylinders.
\end{teorema}

  We stress that both the estimate and the constant in  Theorem \ref{teorema1} are independent of $n$. 

In the second part of the present paper, we focus our study on the case $f\equiv 1$ and we  obtain some quantitative estimates. We define the following functional
\begin{equation}
\label{funzionale}
    \mathcal{F}_p(\Omega)=\frac{T_p(\Omega)P^{q}(\Omega)}{|\Omega|^{q+1}}, \qquad  q=\frac{p}{p-1},
\end{equation}
which is  scaling invariant, since
for every $t>0$
$$\abs{t\Omega}= t^n \abs{\Omega}, \qquad P(t\Omega)= t^{n-1}P(\Omega) $$
and
\begin{equation*}
   T_{p}(t\Omega)=
   t^{n+q}T_{p}(\Omega).
\end{equation*}
We can rewrite inequality \eqref{polyatype}, in the case $f\equiv 1$, as follows
\begin{equation*}
    \mathcal{F}_p(\Omega) \geq c_p.
\end{equation*}
From Theorem \ref{teorema1} follows that  along a sequence of thinning cylinders  $\{\Omega_l\}_{l\in \mathbb{N}}$ defined in \eqref{thin_rect}, we have
\begin{equation*}
    \mathcal{F}_p(\Omega_l) \xrightarrow{l \to 0} c_p . 
\end{equation*}
This leads to the following stability issue: if $\mathcal{F}_p(\Omega)$ is close to $c_p$, can we say that $\Omega$ is close in some sense to a cylinder?
 The following result  gives us information on the nature of the geometry of $\Omega$: when \nobreak{$\mathcal{F}_p(\Omega)-c_p$} is sufficiently small, the set $\Omega$ is a thin domain, in the sense that the ratio $w_\Omega/{\rm diam}(\Omega)$ is small.

The main novelty of the present paper consists indeed in the following quantitative results of the P\'olya estimates \eqref{pol} and the  P\'olya type lower bound in  \eqref{fgl} by means of suitable deficits. 
For completeness, we recall some standard references about isoperimetric quantitative results, see for example \cite{fusco2008sharp, fusco2017stability, brasco_velichkov, brasco2017spectral, gavitone2020quantitative, paoli2020stability}.
The main difference between these results and ours is that the equality in P\'olya's estimates is achieved asymptotically for a sequence of thinning cylinders. Hence, the proof of quantitative result must take into account that we do not have a minimum, as in the classical isoperimetric stability results.

\begin{teorema} \label{teorema2}
Let $\Omega$ be a non-empty, bounded, open and convex set of $\R^n$  and let $f\equiv 1$. Then, 
\begin{equation}\label{quantitative_width}
 \mathcal{F}_p(\Omega)  -c_p \ge K(n,p) \bigg(\frac{w_{\Omega}}{\diam (\Omega)}\bigg)^{n-1},
\end{equation}
where $K(n,p)$ is a positive constant depending only on $p$ and the dimension of the space $n$.
In particular, in the case $n=2$,
the exponent  of the quantity $\displaystyle \frac{w_{\Omega}}{\diam (\Omega)} $ is sharp.
\end{teorema}

We prove a second quantitative result in the case $p=2$ and $n=2$.

\begin{teorema}\label{teorema3}
Let $\Omega$ be a non-empty, bounded, open and convex set in $\R^2$, let $f\equiv 1$ and let $p=2$. Then, there exists a positive constant $\tilde{K}$ such that
\begin{equation}
\label{quantitativa2}
    \mathcal{F}_2(\Omega) -c_2=
    \frac{T(\Omega)P^2(\Omega)}{|\Omega|^3}-\frac{1}{3} \ge \tilde{K} \left(\frac{|\Omega \bigtriangleup Q|}{|\Omega|}\right)^3,
\end{equation}
where $\Omega \bigtriangleup Q$ denotes the symmetric difference between $\Omega$ and a rectangle $Q$ with sides $P(\Omega)/ 2$ and $w_{\Omega}$ containing $\Omega$.
\end{teorema}

We prove the  first quantitative result \eqref{quantitative_width} starting from the proof of Theorem \ref{teorema1}, estimating the remainder term with geometric quantities. 
Finally, in order to prove Theorem \ref{teorema3}, we  use Steiner formulas (see Section \ref{innere} and the references therein).

  The paper is organized as follows. In Section $2$ we recall some basic notions, definitions and we recall some classical results, focusing in particular on the class of convex sets. In Section $3$ we prove Theorem \ref{teorema1} and, finally,  Section $4$ is dedicated to the proof of the quantitative results 
  when $f\equiv 1$.

\section{Notations and Preliminaries}
\label{section_notion}
\subsection{Notations and basic facts} 
 Throughout this article, $|\cdot|$ will denote the Euclidean norm in $\mathbb{R}^n$,
 while $\cdot$ is the standard Euclidean scalar product for  $n\geq2$. By $\mathcal{H}^k(\cdot)$, for $k\in [0,n)$, we denote the $k-$dimensional Hausdorff measure in $\mathbb{R}^n$.
 
 The perimeter of $\Omega$ in $\mathbb{R}^n$ will be denoted by $P(\Omega)$ and, if $P(\Omega)<\infty$, we say that $\Omega$ is a set of finite perimeter. In our case, $\Omega$ is a bounded, open and convex set; this ensures us that $\Omega$ is a set of finite perimeter and that $P(\Omega)=\mathcal{H}^{n-1}(\partial\Omega)$. Moreover, if $\Omega$ is an open set with Lipschitz boundary, it holds

 \begin{teorema}[Coarea formula]
 Let $\Omega \subset \mathbb{R}^n$ be an open set. Let $f\in W^{1,1}_{\text{loc}}(\Omega)$ and let $u:\Omega\to\R$ be a measurable function. Then,
 \begin{equation}
   \label{coarea}
   {\displaystyle \int _{\Omega}u(x)|\nabla f(x)|dx=\int _{\mathbb {R} }dt\int_{\Omega\cap f^{-1}(t)}u(y)\, d\mathcal {H}^{n-1}(y)}.
 \end{equation}
 \end{teorema}
Some references for results relative to the sets of finite perimeter and for the coarea formula are, for instance, \cite{maggi2012sets,ambrosio2000functions}.

 We give now the definition of the support function of a convex set and minimal width (or thickness) of a convex set.
\begin{definizione}\label{support}
  Let $\Omega$ be a bounded, open and convex set of $\mathbb{R}^n$. The support function of $\Omega$ is defined as
  \begin{equation*}
    h_\Omega(y)=\sup_{x\in \Omega}\left(x\cdot y\right), \qquad y\in \mathbb{R}^n .
  \end{equation*}
\end{definizione}

\begin{definizione}
  Let $\Omega$ a bounded, open and convex set of $\mathbb{R}^n$, the width of $\Omega$ in the direction $y \in \mathbb{R}$ is defined as 
  \begin{equation*}
    \omega_{\Omega}(y)=h_{\Omega}(y)+h_{\Omega}(-y)
    \end{equation*}
 and the minimal width of $\Omega$ as
\begin{equation*}
    w_\Omega=\min\{  \omega_{\Omega}(y)\,|\; y\in\mathbb{S}^{n-1}\}.
\end{equation*}
\end{definizione}

We recall the following estimate, which is proved in \cite{fenchel_bonnesen} in the planar case and is generalized in \cite{brasco_2020_principal_frequencies} to all dimensions. 
 \begin{prop}
 Let $\Omega$ be a non-empty bounded, open and convex set of $\mathbb{R}^n$. Then, 
 \begin{equation}\label{convex_estimates}
 \dfrac{1}{n}  \leq\dfrac{|\Omega|}{P(\Omega) R_\Omega}  < 1. 
 \end{equation}
The upper bound is sharp on a sequence of thinning cylinders, while the lower bound is sharp, for example, on balls. Moreover, for $n=2$, any circumscribed polygon, that is a polygon whose incircle touches all the sides, verifies the lower bound with the equality sign. 
 \end{prop}

In the planar case the following inequalities hold true (see as a reference \cite{ scott_convex_2000, scott_family, santalo_sobre}).
\begin{prop}
Let $\Omega$ be a bounded, open and convex set of $\mathbb{R}^2$. Then,

\begin{equation}\label{width_inradius}
    2\leq \frac{w_\Omega}{R_{\Omega}} \leq 3.
\end{equation}
  The upper bound is achieved by equilateral triangles and the lower bound is achieved by disks.\\
  Moreover,
\begin{equation}\label{scott}
 \left(w_{\Omega}-2R_{\Omega} \right) P(\Omega)\leq \dfrac{2}{\sqrt{3}}w^2_{\Omega},
\end{equation}
with equality holding for equilateral triangles, and
\begin{equation}\label{santalo}
|\Omega|\leq R_{\Omega}\left( P(\Omega)-\pi R_{\Omega} \right),
\end{equation}
with equality holding for the stadii (convex hull of two identical disjoint balls).\\
Eventually,
\begin{equation} \label{diam_per}
    2\diam (\Omega)<P(\Omega)\leq \pi \diam(\Omega),
\end{equation}
where the lower bound is asymptotically achieved by a sequence of thinning rectangles and the upper bound by sets of constant width.

\end{prop}

\subsection{Inner parallel sets} \label{innere}
 Let $\Omega$ be a non-empty, bounded, open and convex set of $\R^n$. We defined the distance function from the boundary in \eqref{disti}
 and we will denote it by $d(\cdot)$. We remark that the distance function is concave, as a consequence of the convexity of $\Omega$.

The superlevel sets of the distance function
\begin{equation}\label{inner}
    \Omega_t=\set{x\in \Omega \, : \, d(x)>t}, \qquad t\in[0, R_\Omega]
\end{equation}
 are called \emph{inner parallel sets}, where $R_\Omega$ is the inradius of $\Omega$, and we use the following notations:
 \begin{equation}
 \label{notations}
 \mu(t)= \abs{\Omega_t}, \qquad P(t)=P(\Omega_t)\qquad t\in[0, R_\Omega].
 \end{equation}
 By coarea formula \eqref{coarea}, recalling that $\abs{\nabla d}=1$ almost everywhere, we have 
 $$\mu(t)=\int_{\Set{d>t}} \,  dx = \int_{\Set{d>t}} \frac{\abs{\nabla d}}{\abs{\nabla d}} \,  dx= \int_{t}^{R_\Omega} \frac{1}{\abs{\nabla d}} \int_{\Set{d=s}} d\mathcal{H}^{n-1} \, ds= \int_{t}^{R_\Omega} P(s)\; ds;$$
hence, the function $\mu(t)$ is absolutely continuous, decreasing and its derivative is $\mu'(t)=-P(t)$ almost everywhere. Moreover, it is possible to prove that the perimeter $P(t)$ is non increasing and absolutely continuous, as a consequence of the concavity of the distance function and the Brunn-Minkowski inequality for the perimeter (see \cite{schneider} as a reference).
 
 Finally, let us consider the case $n=2$. For $\Omega$ non-empty bounded, open and convex set of $\mathbb{R}^2$, the Steiner formulas for the inner parallel sets hold (see \cite{Nagy}):
  \begin{equation}\label{steiner1}
     P(t)\leq P(\Omega)-2\pi t \qquad \forall t \in [0,R_\Omega],
 \end{equation}
 \begin{equation}\label{steiner2}
     \mu(t)\geq |\Omega|-P(\Omega) t+\pi t^2 \qquad \forall t \in [0,R_\Omega],
 \end{equation}
 equality holding in both \eqref{steiner1} and \eqref{steiner2} for the stadii (see \cite{fragala_gazzola_lamboley}).
 
 As a consequence of the Alexandrov-Fenchel inequality and the isoperimetric inequality for the quermassintegrals (see \cite{schneider}), we have 
 \begin{equation}\label{dimn}
      -P'(t) \ge  n(n-1)\omega_n^{\frac{1}{n-1}}\bigg(\frac{P(t)}{n}\bigg)^{\frac{n-2}{n-1}}, 
  \end{equation}
that, for $n=2$,  reads
\begin{equation}\label{per_lower}
    -P'(t) \ge 2\pi,
\end{equation}
with  equality if $\Omega$ is a ball or a stadium.

\section[Proof of Theorem 1.1]{Proof of Theorem \ref{teorema1}}

In this Section we prove Theorem \ref{teorema1}.
Since the proof is quite long, we split it in two parts: firstly we prove inequality \eqref{polyatype} and, then, we prove its sharpness.
 \subsubsection*{Step 1: proof of inequality \eqref{polyatype} in Theorem \ref{teorema1}}
\begin{proof}
 Let us choose in the variational characterization \eqref{variational} $\varphi(x)= g(d(x))$ as a test function, where $g$ is a positive and non-decreasing function in $W^{1,p}([0,R_\Omega])$ such that $g(0)=0$. Then,
by coarea formula \eqref{coarea}, 
\begin{equation}
\label{fgP}
    \int_\Omega f(d(x,\partial \Omega))\varphi(x) \, dx = \int_0^{R_\Omega}f(t)g(t) P(t) \, dt
\end{equation}
and
\begin{equation}
\label{g'P}
    \int_\Omega \abs{\nabla \varphi(x)}^p \, dx = \int_0^{R_\Omega}g'^p(t)  P(t) \, dt.
\end{equation}
By \eqref{variational}, \eqref{fgP} and \eqref{g'P} we have
\begin{equation}
\label{tfp_estim}
    T_{f,p}(\Omega) \geq \frac{\displaystyle{\left(\int_0^{R_\Omega}f(t)g(t) P(t) \, dt\right)^\frac{p}{p-1}}}{\left(\displaystyle{\int_0^{R_\Omega}g'^p(t) P(t) \, dt}\right)^{\frac{1}{p-1}}}.
\end{equation}
  Now, if we define the following measure
\begin{equation*}
    \mu_f(E) = \int_E f(d(x)) \, dx,
\end{equation*}
we have
\begin{equation}\label{weig_measure}
    \mu_f(t):=\mu_f(\Omega_t) = \int_{\Omega_t} f(d(x)) \, dx= \int_t^{R_\Omega} f(s) P(s)\, ds.
\end{equation}
Since  $f(s)P(s)$ is a decreasing function, we get 
\begin{equation}
\label{fmisure_vs_inradius}
    \mu_f(t) \leq (R_\Omega - t) f(t) P(t).
\end{equation}
From \eqref{weig_measure}, we have 
\begin{equation}
\label{derivate_volum}
    -\mu'_f(t)= f(t) P(t) \qquad \text{ a.e. } t \in [0, R_\Omega].
\end{equation}
Using \eqref{fgP}, \eqref{derivate_volum} and integrating by parts, we obtain
\begin{equation*}
   \int_0^{R_\Omega}f(t)g(t) P(t) \, dt=-\int_0^{R_\Omega}g(t)\mu'_f(t) \, dt=\int_0^{R_\Omega}g'(t)\mu_f(t)\, dt.
\end{equation*}
Consequently, \eqref{tfp_estim} becomes
\begin{equation*}
    T_{f,p}(\Omega) \geq \frac{\displaystyle{\left(\int_0^{R_\Omega}g'(t)\mu_f(t)\, dt\right)^{\frac{p}{p-1}}}}{\left(\displaystyle{\int_0^{R_\Omega}g'^p(t) P(t) \, dt}\right)^{\frac{1}{p-1}}}.
\end{equation*}
We can choose 
\begin{equation*}
    g(t)=\int_0^t\left( \dfrac{\mu_f(s)}{P(s)}\right)^{1/(p-1)}ds
    \end{equation*}
and we observe that $g\in W^{1,p}([0, R_{\Omega}])$, since, using \eqref{fmisure_vs_inradius}, we have
\begin{gather*}
   g(t) \leq \int_0^{R_\Omega} (R_\Omega -s)^{\frac{1}{p-1}} f(s)^{\frac{1}{p-1}}\,ds \leq \norma{f}_{L^\infty}^{\frac{1}{p-1}} R_\Omega^{\frac{p}{p-1}}\in L^p([0, R_{\Omega}]),
    \\
    g'(t) \leq \norma{f}_{L^\infty}^{\frac{1}{p-1}} R_\Omega^{\frac{1}{p-1}} \in L^p([0, R_{\Omega}]).
\end{gather*}
So, we have
\begin{equation}\label{lower_tor1}
    T_{f,p}(\Omega) \geq
\int_0^{R_\Omega}\frac{\mu_f^{\frac{p}{p-1}}(t)}{P^{\frac{1}{p-1}}(t)}\, dt =-\frac{p-1}{2p-1}\int_0^{R_\Omega}\frac{(\mu^{\frac{2p-1}{p-1}}_f(t))'}{f(t)P^{\frac{p}{p-1}}(t)}\, dt.
\end{equation}
Let us set $c_p =(p-1)/(2p-1)$. Since $f(s)$ is a non-negative and non-increasing function, integrating by parts in \eqref{lower_tor1}, we get
\begin{equation} \label{lower-tor2}
\begin{aligned}
    T_{f,p}(\Omega)&\ge -c_p\int_0^{R_\Omega}\frac{(\mu^{\frac{2p-1}{p-1}}_f(t))'}{f(t)P^{\frac{p}{p-1}}(t)}\, dt =-c_p\frac{\mu^{\frac{2p-1}{p-1}}_f(t)}{f(t)P^{\frac{p}{p-1}}(t)} \Bigg\rvert_0^{R_\Omega} + \\
    &-c_p\int_0^{R_\Omega}\frac{\mu^{\frac{2p-1}{p-1}}_f(t)}{f^2(t)P^{\frac{2p}{p-1}}(t)} \left( f'(t)P^{\frac{p}{p-1}}(t)+ \frac{p}{p-1} f(t)P^{\frac{1}{p-1}}(t)P'(t)\right)\, dt\\
& \geq c_p\frac{\mu^{\frac{2p-1}{p-1}}_f(\Omega)}{f(0)P^{\frac{p}{p-1}}(\Omega)} + \frac{c_p}{ P^{\frac{p}{p-1}}(\Omega)} \int_0^{R_\Omega}\frac{\mu^{\frac{2p-1}{p-1}}_f(t)}{f^2(t)} (-f'(t))\, dt,
\end{aligned}
\end{equation}
where in the  last inequality we use \eqref{fmisure_vs_inradius} and the fact that $P'(t)\le 0$.
Now, since $f(s)$ is non-increasing, we obtain the desired estimate 

\begin{equation}\label{lower-tor3}
    T_{f,p}(\Omega)\ge c_p\frac{\mu^{\frac{2p-1}{p-1}}_f(\Omega)}{f(0)P^{\frac{p}{p-1}}(\Omega)}. 
\end{equation}

\end{proof}

\subsubsection*{Step 2: proof of the sharpness of \eqref{polyatype} }
\begin{proof}
We prove that inequality \eqref{polyatype} is sharp and that the optimum is asymptotically achieved by the sequence of thinning cylinders $\Omega_l$ with unitary measure, as defined in \eqref{thin_rect}, that is
    \begin{equation*}
    \Omega_l = l^{-\frac{1}{n-1}}C \times \left(-\frac{l}{2}, \frac{l}{2}\right)
\end{equation*}
where $C\subseteq \R^{n-1}$ is a bounded, open and convex set with unitary $(n-1)-$measure.
It is easy to verify that, for $n\geq 3$,
\begin{equation}
   \begin{aligned}
\label{perimeter_cylinder}
    P(\Omega_l)&= 2 \mathcal{H}^{n-1} (l^{-\frac{1}{n-1}} C) + l \mathcal{H}^{n-2} (\partial( l^{-\frac{1}{n-1}}C) )\\
    &=2l^{-1}+ l^{\frac{1}{n-1}} \mathcal{H}^{n-2}(\partial C),
    \end{aligned} 
\end{equation}
and we observe that, in the case $n=2$, we have that $\mathcal{H}^{n-2}(\partial C)=2$.

 Let $u$ be the solution to the following $p$-torsion problem
\begin{equation*}
\begin{cases}
-\Delta_p u = 1 &\text{in } \Omega_l\\
u=0 &\text{on } \partial \Omega_l, 
\end{cases}
\end{equation*}
such that 
\begin{equation*}
    \int_{\Omega_l} u\;dx=T_p(\Omega_l),
\end{equation*}
and let us consider the following function, depending only on the last component $x_n$ of $x\in \R^n$,
\begin{equation*}
    v(x) = \frac{p-1}{p} \left[ \left(\frac{l}{2}\right)^{\frac{p}{p-1}} - \abs{x_n}^{\frac{p}{p-1}}\right],
\end{equation*}
satisfying
\begin{equation*}
\begin{cases}
-\Delta_p v = 1 &\text{in } \Omega_l\\
v\ge0 &\text{on } \partial \Omega_l.
\end{cases}
\end{equation*}
The comparison principle, see \cite{lindqvist2017notes}, ensures that $u\le v$ in $\Omega_l$ and, as a consequence, 
\begin{equation}
\label{p-tortionup}
\begin{aligned}
T_{p}(\Omega_l) &= \int_{\Omega_l} u\,dx \le \int_{\Omega_l}v\,dx  = 
\\&=\frac{p-1}{p} \int_{l^{-\frac{1}{n-1} }C} \int_{-\frac{l}{2}}^{\frac{l}{2}}\left[ \left(\frac{l}{2}\right)^{\frac{p}{p-1}} - \abs{x_n}^{\frac{p}{p-1}}\right] \,dx_n\,d\mathcal{H}^{n-1} \\
&=2\; \frac{p-1}{p} l^{-1} \int_0^{\frac{l}{2}} \left[ \left(\frac{l}{2}\right)^{\frac{p}{p-1}} - x_n^{\frac{p}{p-1}}\right]\,dx_n\\
&=2\; \frac{p-1}{p} \left[ 1 - \frac{p-1}{2p-1}\right]l^{-1} \left(\frac{l}{2}\right)^{\frac{2p-1}{p-1}} = 2c_p l^{-1} \left(\frac{l}{2}\right)^{\frac{2p-1}{p-1}}.
\end{aligned}
\end{equation}
By \eqref{p-tortionup} and \eqref{perimeter_cylinder}, we have 
\begin{equation*}
    T_p(\Omega_l) P^{\frac{p}{p-1}}(\Omega_l) \le 
    2c_p l^{-1} \left(\frac{l}{2}\right)^{\frac{2p-1}{p-1}} \left( 2l^{-1}+ l^{\frac{1}{n-1}} \mathcal{H}^{n-2}(\partial C)\right)^{\frac{p}{p-1}} = c_p \left( 1+ \frac{l^{\frac{n}{n-1}}}{2}\mathcal{H}^{n-2}(\partial C)\right)^{\frac{p}{p-1}}.
\end{equation*}
Now, since $f(x)\leq f(0)$, we have that, for every bounded, open and convex set $\Omega$,
\begin{equation}\label{fTorsion}
    T_{f,p}(\Omega)\leq  f^{\frac{p}{p-1}}(0) T_p(\Omega).
\end{equation}
It follows that 
\begin{equation}
\label{p-tortion.per}
\begin{aligned}
    T_{f,p}(\Omega_l)P^{\frac{p}{p-1}}(\Omega_l)&\le f^{\frac{p}{p-1}}(0) T_p(\Omega_l) P^{\frac{p}{p-1}}(\Omega_l)\\
     &\le c_p f^{\frac{p}{p-1}}(0) \left( 1+ \frac{l^{\frac{n}{n-1}}}{2}\mathcal{H}^{n-2}(\partial C)\right)^{\frac{p}{p-1}}.
\end{aligned}
\end{equation}
Moreover we observe that, if $f$ never vanishes, we can use its monotonicity property to bound $\mu_f$ from below in the following way:

\begin{equation*}
    \mu_f(\Omega) = \int_{\Omega}f(d(x))\,dx \ge f(R_{\Omega})|\Omega|,
\end{equation*}
obtaining 
\begin{equation}\label{new}
    T_{f,p}(\Omega)\ge c_p\frac{f^{\frac{2p-1}{p-1}}(R_{\Omega})|\Omega|^{\frac{2p-1}{p-1}}}{f(0)P^{\frac{p}{p-1}}(\Omega)}.
\end{equation}
Joining \eqref{new} and \eqref{p-tortion.per}, we obtain
\begin{equation*}
    c_p\frac{f^{\frac{2p-1}{p-1}}(R_{\Omega_l})}{f(0)}\le T_{f,p}(\Omega_l)P^{\frac{p}{p-1}}(\Omega_l)\le c_p f^{\frac{p}{p-1}}(0) \left( 1+ \frac{l^{\frac{n}{n-1}}}{2}\mathcal{H}^{n-2}(\partial C)\right)^{\frac{p}{p-1}}.
\end{equation*}
Eventually, passing to the limit when $l\to 0$, observing that $\displaystyle\lim_{l\to 0} R_{\Omega_l}=0$ and that $f$ is continuous, we have
\begin{equation*}
    T_{f,p}(\Omega_l)P^{\frac{p}{p-1}}(\Omega_l) \longrightarrow c_p f^{\frac{p}{p-1}}(0).
\end{equation*}
\end{proof}

\begin{oss}
If we assume that $f:[0,R_\Omega]\to [0,+\infty[$ is a function in $L^{\infty}([0,R_\Omega])$, then, using the variational characterization \eqref{variational} and the result \eqref{gavit} proved in \cite{gavitone_2014}, we have 
\begin{equation}\label{dell}
   T_{f,p}(\Omega)\geq \left(\inf_{t \in [0, R_\Omega]} f(t)\right)^{\frac{p}{p-1}} T_p(\Omega)\geq  \left(\inf_{t \in [0, R_\Omega]} f(t)\right)^{\frac{p}{p-1}} c_p \dfrac{|\Omega|^{\frac{2p-1}{p-1}}}{P(\Omega)^{\frac{p}{p-1}}}
   \end{equation}
and the sharpness of \eqref{dell} can be proved in an analogous way as in \eqref{polyatype}.
\end{oss}

\section{The quantitative results}

\subsubsection*{Proof of Theorem \ref{teorema2}}
\begin{proof}    

Let us start by proving \eqref{quantitative_width} in the case $n=2$. If $f\equiv 1$, \eqref{lower-tor2} becomes
\begin{equation} \label{cc}
    T_p(\Omega) \ge c_p \frac{|\Omega|^{\frac{2p-1}{p-1}}}{P^{\frac{p}{p-1}}(\Omega)}+ c_p \frac{p}{p-1}\int_0^{R_{\Omega}} \left(\frac{\mu(t)}{P(t)}\right)^{\frac{2p-1}{p-1}}(-P'(t))\,dt.
\end{equation}
Joining \eqref{convex_estimates}, \eqref{width_inradius}, \eqref{per_lower} and \eqref{cc}, we have that 
\begin{equation*}
\begin{aligned}
\frac{T_p(\Omega)P^{\frac{p}{p-1}}(\Omega)}{|\Omega|^{\frac{2p-1}{p-1}}}
    -c_p  & > c_p \frac{p}{p-1}\frac{P^{\frac{p}{p-1}}(\Omega)}{|\Omega|^{\frac{2p-1}{p-1}}}\int_0^{R_{\Omega}} \left(\frac{\mu(t)}{P(t)}\right)^{\frac{2p-1}{p-1}}(-P'(t))\,dt \\
    &\ge \frac{\pi}{2^{\frac{p}{p-1}}} \frac{p}{2p-1}  \frac{P^{\frac{p}{p-1}}(\Omega)}{|\Omega|^{\frac{2p-1}{p-1}}} \int_0^{R_{\Omega}} \left(R_\Omega-t\right)^{\frac{2p-1}{p-1}}\,dt \\
    &\ge \frac{\pi}{2^{\frac{p}{p-1}}} \frac{(p-1)p}{(3p-2)(2p-1)}
    \frac{R_\Omega}{P(\Omega)}\left(\frac{R_\Omega P(\Omega)}{|\Omega|} \right)^{{\frac{2p-1}{p-1}}} \\
    &\ge \frac{\pi}{2^{\frac{p}{p-1}}} \frac{(p-1)p}{(3p-2)(2p-1)}
    \frac{R_\Omega}{P(\Omega)}.
    \end{aligned}
\end{equation*}
Hence, by applying \eqref{diam_per} and \eqref{width_inradius} we get
\begin{equation} 
\label{4c}
\mathcal{F}_p(\Omega) -c_p \ge K(2,p)
    \frac{w_\Omega}{\diam(\Omega)},
\end{equation}
where
\begin{equation}
\label{consty}
    K(2,p)=\frac{(p-1)p}{ 2^{\frac{p}{p-1}}3(3p-2)(2p-1)} .
\end{equation}

We now prove that the exponent of  $\displaystyle \frac{w_\Omega}{\diam(\Omega)}$
 in \eqref{4c} is sharp.
In order to do that, 
we only need to find a sequence $\{\Omega_l\}_{l\in \mathbb{N}}$ of convex sets with fixed measure such that
\begin{equation*}
M \frac{w_{\Omega_l}}{\diam(\Omega_l)} \ge \mathcal{F}_p(\Omega_l) -c_p ,
\end{equation*}
for some positive constant $M$.
Let $0<l<1$, we consider the following rectangle
\begin{equation*}
   \Omega_l = \left( -\frac{1}{2l},\frac{1}{2l}\right) \times \left( -\frac{l}{2}, \frac{l}{2}\right)
\end{equation*}
and we notice that its inradius and area are $R_{\Omega_l}= \frac{l}{2}$ and $\abs{\Omega_l}=1$.
Let   $u$ be the unique solution to
\begin{equation*}
    \begin{cases}
    -\Delta_p u = 1 \quad \text{in } \Omega_l \\
    u=0 \qquad \text{on } \partial\Omega_l
    \end{cases}
\end{equation*}
and let us consider the following function
\begin{equation*}
    v(y) = \frac{p-1}{p}\left[ \left(\frac{l}{2}\right)^{\frac{p}{p-1}} -\abs{y}^{\frac{p}{p-1}}\right],
\end{equation*}
which solves
\begin{equation*}
\begin{cases}
-\Delta v = 1 &\text{in } \Omega_l\\
v\ge0 &\text{on } \partial \Omega_l.
\end{cases}
\end{equation*}
The comparison principle gives $u\le v$ in $\Omega_l$ and
\begin{equation*}
    T_p(\Omega_l) = \int_{\Omega_l}u_p\,dx \le \int_{\Omega_l}v\,dx.
 \end{equation*}
 Arguing as in \eqref{p-tortionup}, we have
\begin{equation*}
    \int_{\Omega_l}v\,dx  =  c_p\left(\frac{l}{2}\right)^{\frac{p}{p-1}}.
\end{equation*}
On the other hand, the perimeter of the rectangle is given by
\begin{equation*}
    P(\Omega_l) = \frac{2}{l}\left( 1+l^2\right)
\end{equation*}
and its Taylor expansion with respect to $l>0$ is
\begin{equation*}
    P^{{\frac{p}{p-1}}}(\Omega_l) = \left(\frac{2}{l}\right)^{{\frac{p}{p-1}}}\left( 1+l^2\right)^{{\frac{p}{p-1}}} = \left(\frac{2}{l}\right)^{{\frac{p}{p-1}}}\left( 1+\frac{p}{p-1}l^2 + o(l^2)\right).
\end{equation*}
Using \eqref{width_inradius} and \eqref{diam_per}, we get 

\begin{equation*}
\begin{aligned}
   T_p(\Omega_l)P^{\frac{p}{p-1}}(\Omega_l) - c_p& \le  c_p \left(\frac{l}{2}\right)^{\frac{p}{p-1}} \left(\frac{2}{l}\right)^{{\frac{p}{p-1}}}\left( 1+\frac{p}{p-1}l^2+ o(l^2)\right) - c_p\\
&\le 2c_p \frac{p}{p-1} l^2 \leq16 c_p \frac{ p}{p-1} \frac{R_{\Omega_l}}{P(\Omega_l)}\\
&\leq4 c_p \frac{ p}{p-1} \frac{w_{\Omega_l}}{\diam(\Omega_l)}
\end{aligned}
\end{equation*}
and this concludes the proof in dimension $n=2$.

Let us now prove \eqref{quantitative_width} in the case $n>2$.
If we choose $f\equiv 1$, \eqref{lower-tor2} becomes
\begin{equation} \label{dimn1}
    T_p(\Omega) \ge c_p \frac{|\Omega|^{\frac{2p-1}{p-1}}}{P^{\frac{p}{p-1}}(\Omega)}+ c_p \frac{p}{p-1}\int_0^{R_{\Omega}} \left(\frac{\mu(t)}{P(t)}\right)^{\frac{2p-1}{p-1}}(-P'(t))\,dt.
\end{equation}
Hence, combining \eqref{dimn} and \eqref{dimn1}, we have 
\begin{equation}\label{n}
\frac{T_p(\Omega)P^{\frac{p}{p-1}}(\Omega)}{|\Omega|^{\frac{2p-1}{p-1}}}
    -c_p   \ge k(n,p) \frac{P^{\frac{p}{p-1}}(\Omega)}{|\Omega|^{\frac{2p-1}{p-1}}}\int_0^{R_{\Omega}} \left(\frac{\mu(t)}{P(t)}\right)^{\frac{2p-1}{p-1}}P(t)^{\frac{n-2}{n-1}}\,dt.
\end{equation}
Moreover, from \eqref{convex_estimates}, we obtain that 
\begin{equation}\label{inry}
    P(t)\ge n\omega_n (R_{\Omega}-t)^{n-1}, 
\end{equation}
 and so, using \eqref{inry} in \eqref{n}, we get
\begin{equation}\label{nn}
\begin{split}
    \frac{T_p(\Omega)P^{\frac{p}{p-1}}(\Omega)}{|\Omega|^{\frac{2p-1}{p-1}}}
    -c_p   &\ge  k(n,p) \frac{P^{\frac{p}{p-1}}(\Omega)}{|\Omega|^{\frac{2p-1}{p-1}}}\int_0^{R_{\Omega}} \left(R_{\Omega}-t\right)^{\frac{2p-1}{p-1}+n-2}\,dt \\
    &=k(n,p) \bigg(\frac{R_{\Omega}P(\Omega)}{|\Omega|}\bigg)^{\frac{2p-1}{p-1}}\frac{R_{\Omega}^{n-1}}{P(\Omega)}.
 \end{split}   
\end{equation}
If we combine \eqref{nn} with \eqref{convex_estimates}, with the following estimate (that can be found in \cite{fenchel_bonnesen}):
\begin{equation*}
    \displaystyle R_{\Omega}\ge\begin{cases}
    w_{\Omega} \displaystyle{\frac{\sqrt{n+2}}{2n+2}} & n \,\, \text{even}\\ \\
    w_{\Omega} \displaystyle{\frac{1}{2\sqrt{n}}} & n \,\, \text{odd},
    \end{cases}
\end{equation*}
and with 
$$ \displaystyle{P(\Omega)\le n\omega_n \left(\frac{n}{2n+2}\right)^{\frac{n-1}{2}} \diam(\Omega)^{n-1}},$$ we finally get
\begin{equation*}
    \frac{T_p(\Omega)P^{\frac{p}{p-1}}(\Omega)}{|\Omega|^{\frac{2p-1}{p-1}}}
    -c_p \ge K(n,p) \bigg(\frac{w_{\Omega}}{\diam (\Omega)}\bigg)^{n-1}.
\end{equation*}

\end{proof}


\begin{oss}\label{remarkn}

As far as the sharpness of \eqref{quantitative_width} in the case $n>2$, we conjecture that the sharp exponent is $1$ as in the planar case. Indeed, the minimizing sequence $\{\Omega_l\}$ satisfies
$$T_p(\Omega_l)P^{\frac{p}{p-1}}(\Omega_l)-c_p \approx C \frac{w_{\Omega_l}}{\diam(\Omega_l)}.$$
\end{oss}

\begin{oss}
 As already remarked in the Introduction,  inequality \eqref{quantitative_width} gives information on the set $\Omega$. Indeed, if  
\begin{equation*} 
\mathcal{F}_p(\Omega) -c_p 
\end{equation*}
is small, then the ratio between $w_\Omega$ and $\diam(\Omega) $ has to be necessarily small, i.e. $\Omega$ must be a thin domain.
Moreover, inequality \eqref{quantitative_width} tells us also that the infimum of $\mathcal{F}_p(\Omega)$ is not achieved among bounded, open and convex sets. Assuming by contradiction that there exists a bounded, open and convex set $\tilde{\Omega}$ such that 
\begin{equation*}
\mathcal{F}_p(\tilde{\Omega}) = c_p, 
\end{equation*}
 we have that
\begin{equation*}
    \frac{w_{\tilde{\Omega}}}{\diam({\tilde{\Omega}})} < \varepsilon \,\,\,\,\,\,\,\,\,\,\, \forall \varepsilon >0,
\end{equation*}
which is impossible.
\end{oss}

Theorem \ref{teorema2} only tells us that any minimizing sequence of $\mathcal{F}_p(\cdot)$ is a sequence of thinning domains. On the other hand, Theorem \ref{teorema3} gives us more precise information on the geometry of such minimizing sequence in the planar case.

\subsubsection*{Proof of Theorem \ref{teorema3}}

\begin{proof}
Let $\Omega$ be a non-empty, bounded, open and convex set in $\mathbb{R}^2$ and let us consider a rectangle $Q$ with sides $P(\Omega)/2$ and $w_{\Omega}$ containing $\Omega$. Such a rectangle exists, since it is enough to choose the shorter side of $Q$ parallel to the direction of $w_{\Omega}$
and to recall the lower bound in \eqref{diam_per} (see Figure \ref{fig:M3}).

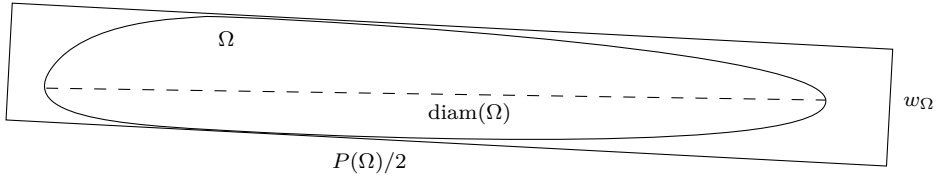
\begin{figure}[h]
\begin{center}
  \begin{tikzpicture}[x=0.75pt,y=0.75pt,yscale=-1,xscale=1]

\draw  (170.5,300) .. controls (187.5,274) and (237.5,274) .. (247.5,273) .. controls (257.5,272) and (554.5,287) .. (557.5,315) .. controls (560.5,343) and (325.5,334) .. (271,331) .. controls (216.5,328) and (153.5,326) .. (170.5,300) -- cycle ;
\draw  [dash pattern={on 4.5pt off 4.5pt}]  (168.5,309) -- (557.5,315) ;
\draw   (151.6,266.24) -- (590.89,289.42) -- (587.79,348.18) -- (148.5,325) -- cycle ;

\draw (356,315) node [anchor=north west][inner sep=0.75pt]   [align=left] {\footnotesize$ \diam(\Omega)$};
\draw (595,312) node [anchor=north west][inner sep=0.75pt]   [align=left] {\footnotesize$w_\Omega$};
\draw (253,279) node [anchor=north west][inner sep=0.75pt]   [align=left] {\footnotesize$\Omega$};
\draw (310,340) node [anchor=north west][inner sep=0.75pt]   [align=left] {\footnotesize$P(\Omega)/2$};
\end{tikzpicture}

\end{center}
\caption{Rectangle with sides $P(\Omega)/2$ and $w_{\Omega}$ containing $\Omega$.} \label{fig:M3}
\end{figure}

Now, let $\sigma>0$ be such that
\begin{gather}
\label{referee1}
   \frac{1}{ 4^3\cdot6} -  \frac{\pi^2}{ 2^3\cdot 3^3}\frac{\sigma^2}{K^2(2)}\geq 0; \\
   \label{referee2}
   \frac{1}{ 3^3\cdot 6} -\frac{\pi}{48}\frac{\sigma}{K(2))} -\frac{\pi^2}{ 2^5\cdot 3}\frac{\sigma^2}{K^2(2)}\geq 0;\\
   \label{referee3}
   \frac{\pi}{4}-\frac{\pi }{2\sqrt{3}}\frac{\sigma}{K(2)}\geq \frac{4}{3 \sqrt{3}},
\end{gather}
where $K(2):=K(2,2)$ is the constant defined in \eqref{consty}. If
$$
\frac{T(\Omega)P^2(\Omega)}{|\Omega|^3}-\frac{1}{3} \geq\sigma,
$$
 then, by \eqref{width_inradius} and \eqref{convex_estimates}, we have 
$$
 \frac{\abs{Q\bigtriangleup \Omega}}{\abs{\Omega}}=\left(\frac{P(\Omega) w_\Omega}{2\abs{\Omega}}-1\right)  \leq \left(\frac{3}{2}\frac{P(\Omega) R_\Omega}{\abs{\Omega}}-1\right)\leq 2.
$$
So, it follows that
\begin{equation*}
\frac{T(\Omega)P^2(\Omega)}{|\Omega|^3}-\frac{1}{3} \geq \frac{\sigma}{2^3}2^3\geq \frac{\sigma}{2^3}\left(\frac{\abs{Q\bigtriangleup \Omega}}{|\Omega|} \right)^3.
\end{equation*}
On the other hand, let us assume that
\begin{equation}
\label{<sig}
\frac{T(\Omega)P^2(\Omega)}{|\Omega|^3}-\frac{1}{3} <\sigma.
\end{equation}
By Theorem \ref{teorema2}, we have that 
    \begin{equation}
    \label{fromteorem1}
     \frac{w_\Omega}{\diam (\Omega)} \leq \frac{1}{K(2)}\left[\frac{T(\Omega)P^2(\Omega)}{|\Omega|^3}-\frac{1}{3} \right]< \frac{\sigma}{K(2)},
    \end{equation}
    and we observe that, by the choice of $\sigma$ made in \eqref{referee1}-\eqref{referee3}, a ball cannot satisfy \eqref{<sig}.
    
   Now,  arguing as in \eqref{lower_tor1} with $f\equiv 1$ and $p=2$, we know that
\begin{equation}
\label{polyas}
    T(\Omega) \ge \int_0^{R_\Omega}\frac{\mu^2(t)}{P(t)}\,dt.
\end{equation}
We set $\displaystyle{\rho=\dfrac{P^2(\Omega)}{4 \pi}  - \abs{\Omega}}$ and $p_R=P(\Omega)- 2 \pi R_\Omega$ and we observe that they are both strictly positive by the isoperimetric inequality and the monotonicity of the perimeter, respectively.
Using inequalities \eqref{steiner1} and \eqref{steiner2} in \eqref{polyas}, we have that 
\begin{equation}
\label{tp^2}
\begin{split}
    T(\Omega)P^2(\Omega) &\ge P^2(\Omega) \int_0^{R_\Omega} \frac{(|\Omega|-P(\Omega)t+\pi t^2)^2}{P(\Omega)-2\pi t}\,dt\\
    &= P^2(\Omega) \int_0^{R_\Omega} \frac{1}{P(\Omega) -2 \pi t}\left(\frac{\left(P(\Omega) -2 \pi t\right)^2}{4\pi} -\left(\frac{P^2(\Omega)}{4 \pi} - \abs{\Omega}\right)\right)^2 \,dt\\
    &= P^2(\Omega) \int_0^{R_\Omega} \left(\frac{\left(P(\Omega) -2 \pi t\right)^3}{(4\pi)^2} - \frac{\rho}{2 \pi } \left(P(\Omega) -2 \pi t\right) +\frac{\rho^2}{P(\Omega) -2 \pi t}\right)\, dt\\
   &= \frac{P^2(\Omega) }{2\pi} \left(\frac{P^4(\Omega)-p_R^4}{4(4\pi)^2} - \frac{\rho}{4 \pi } \left(P^2(\Omega) -p_R^2\right) -\rho^2 \log\left(1-\frac{2\pi R_\Omega}{P(\Omega)}\right)\right),
\end{split}
\end{equation}
 and, using Newton's formula and the Taylor series for the logarithm, we get
\begin{equation}
\label{taylor}
\begin{aligned}
    &P^2(\Omega)-p^2_R= 4\pi R_\Omega P(\Omega)-4\pi^2 R^2_\Omega;\\
   & P^4(\Omega)-p^4_R= 8\pi R_\Omega P^3(\Omega) -24 \pi^2 R_\Omega^2 P^2(\Omega)+32 \pi^3 R^3_\Omega P(\Omega )- 16 \pi^4 R^4_\Omega;\\
   & -\log\left(1-\frac{2\pi R_\Omega}{P(\Omega)}\right)=\sum_{i=1}^{\infty}\frac{1}{i} \left(\frac{2\pi R_\Omega}{P(\Omega)}\right)^i \geq\frac{2\pi R_\Omega}{P(\Omega)}+ \frac{2\pi^2 R^2_\Omega}{P^2(\Omega)} + \frac{8}{3}\frac{\pi^3 R^3_\Omega}{P^3(\Omega)} + \frac{4\pi^4 R^4_\Omega}{P^4(\Omega)}.
\end{aligned}
\end{equation}
 By \eqref{taylor} and \eqref{tp^2}, dividing by $\abs{\Omega}^3$ and subtracting $1/3$, we have
\begin{equation} \label{badineq}
\begin{aligned}
     \frac{T(\Omega)P^2(\Omega)}{\abs{\Omega}^3} - \frac{1}{3} \ge & \frac{1}{3}\bigg(\frac{P(\Omega)R_\Omega}{\abs{\Omega}}-1\bigg)^3 +\pi \frac{R_\Omega^2}{\abs{\Omega}^2}  \left(\abs{\Omega}-\frac{2}{3} P(\Omega)R_\Omega \right)\\
     &+\frac{4}{3}\pi^2 \frac{R_\Omega^3}{P(\Omega)\abs{\Omega}^2} \left(\abs{\Omega}-\frac{3}{4} P(\Omega)R_\Omega \right).
\end{aligned}
\end{equation}
As an intermediate step we want to prove the following inequality: 
\begin{equation}
\label{1/6 prw}
\begin{multlined}
    \frac{1}{3}\bigg(\frac{P(\Omega)R_\Omega}{\abs{\Omega}}-1\bigg)^3 +\pi \frac{R_\Omega^2}{\abs{\Omega}^2}  \left(\abs{\Omega}-\frac{2}{3} P(\Omega)R_\Omega \right)\\
    +\frac{4}{3}\pi^2 \frac{R_\Omega^3}{P(\Omega)\abs{\Omega}^2} \left(\abs{\Omega}-\frac{3}{4} P(\Omega)R_\Omega \right) \geq \frac{1}{6}\bigg(\frac{P(\Omega)R_\Omega}{\abs{\Omega}}-1\bigg)^3,
\end{multlined}
\end{equation}
that, combined with \eqref{badineq}, implies
\begin{equation}
\label{quantitative_R}
   \frac{T(\Omega)P^2(\Omega)}{|\Omega|^3}-\frac{1}{3} \ge \frac{1}{6} \left(\frac{P(\Omega)R_\Omega}{|\Omega|}-1 \right)^3,
\end{equation}
where we choose the constant $1/6$ as an arbitrary constant less then $1/3$.
In particular, \eqref{1/6 prw} is equivalent to 
\begin{equation}
\label{no_taylor}
    \frac{1}{6}\big(P(\Omega)R_\Omega-|\Omega|\big)^3 +\pi R_\Omega^2 \abs{\Omega} \left(\abs{\Omega}-\frac{2}{3} P(\Omega)R_\Omega \right)+\frac{4}{3}\pi^2 \frac{R_\Omega^3}{P(\Omega)}\abs{\Omega} \left(\abs{\Omega}-\frac{3}{4} P(\Omega)R_\Omega \right) \geq 0.
\end{equation}
 In order to prove \eqref{no_taylor}, we distinguish three cases:
\begin{enumerate}
    \item[1)]  if $\displaystyle{\abs{\Omega}\geq \frac{3}{4} P(\Omega)R_\Omega }$, then \eqref{no_taylor} is trivial, since the left hand side is the sum of positive quantities;
    \item[2)]  if $\displaystyle{\frac{2}{3}P(\Omega)R_\Omega \leq\abs{\Omega}< \frac{3}{4} P(\Omega)R_\Omega}$, using \eqref{width_inradius}, \eqref{diam_per}, \eqref{referee1} and \eqref{fromteorem1}, we have
    \begin{equation}
        \label{2/3 and 3/4}
    \begin{split}
       \frac{1}{6}\big(P(\Omega)R_\Omega-|\Omega|\big)^3 &+\pi R_\Omega^2 \abs{\Omega} \left(\abs{\Omega}-\frac{2}{3} P(\Omega)R_\Omega \right) +\frac{4}{3}\pi^2 \frac{R_\Omega^3}{P(\Omega)}\abs{\Omega} \left(\abs{\Omega}-\frac{3}{4} P(\Omega)R_\Omega \right) \\ &\geq P^3(\Omega)R^3_\Omega\left(\frac{1}{ 4^3\cdot6} -\frac{2\pi^2}{3^3}\frac{R^2_\Omega}{P^2(\Omega)}\right)\\ &\geq
       P^3(\Omega)R^3_\Omega\left(\frac{1}{ 4^3\cdot6} -
        \frac{\pi^2}{ 2^3\cdot 3^3}\frac{w^2_\Omega}{\diam^2(\Omega)}\right)
         \\
        &\geq  P^3(\Omega)R^3_\Omega\left(\frac{1}{ 4^3\cdot6} -
        \frac{\pi^2}{ 2^3\cdot 3^3}\frac{\sigma^2}{K^2(2)}\right) \geq 0.
    \end{split}
    \end{equation}  
     \item[3)] if $\displaystyle{\frac{1}{2}P(\Omega)R_\Omega \leq\abs{\Omega}< \frac{2}{3} P(\Omega)R_\Omega}$, arguing as before, we have
     \begin{equation}
        \label{less2/3}
    \begin{split}
       \frac{1}{6}\big(P(\Omega)R_\Omega-|\Omega|\big)^3 &+\pi R_\Omega^2 \abs{\Omega} \left(\abs{\Omega}-\frac{2}{3} P(\Omega)R_\Omega \right) +\frac{4}{3}\pi^2 \frac{R_\Omega^3}{P(\Omega)}\abs{\Omega} \left(\abs{\Omega}-\frac{3}{4} P(\Omega)R_\Omega \right) \\ &\geq P^3(\Omega)R^3_\Omega\left(\frac{1}{ 3^3\cdot 6} -\frac{\pi}{48}\frac{w_\Omega}{\diam(\Omega)} -\frac{\pi^2}{ 2^5\cdot 3}\frac{w^2_\Omega}{\diam^2(\Omega)}\right)\\
       &\geq P^3(\Omega)R^3_\Omega\left(\frac{1}{ 3^3\cdot 6} -\frac{\pi}{48}\frac{\sigma}{K(2)} -\frac{\pi^2}{ 2^5\cdot 3}\frac{\sigma^2}{K^2(2)}\right)\geq 0.
    \end{split}
    \end{equation} 
\end{enumerate}
So, we have proved the intermediate step \eqref{quantitative_R}.
Now, by combining \eqref{quantitative_R} and \eqref{scott}, we deduce
\begin{equation}
\label{pr_vs_pw}
       \frac{T(\Omega)P^2(\Omega)}{|\Omega|^3}-\frac{1}{3} \ge\frac{1}{6} \left[\frac{P(\Omega)R_\Omega}{|\Omega|}-1 \right]^3\ge\frac{1}{6} \left[\frac{P(\Omega) w_\Omega}{2|\Omega|}-1 -\frac{1}{\sqrt{3}}\frac{w_\Omega^2}{\abs{\Omega}} \right]^3.
\end{equation}
Using \eqref{santalo}, \eqref{scott}, \eqref{fromteorem1}  and \eqref{referee3}, we have 
\begin{equation} \label{pw-m}
\begin{split}
    \frac{P(\Omega) w_\Omega}{2\abs{\Omega}}-1&\ge
    \frac{P(\Omega) R_\Omega}{\abs{\Omega}}-1\ge \pi \frac{R_\Omega^2}{\abs{\Omega}}\ge \frac{\pi}{\abs{\Omega}}\left(\frac{w_\Omega}{2}- \frac{w_\Omega^2}{\sqrt{3}P(\Omega)}\right)^2 \\
    &=\frac{w_\Omega^2}{\abs{\Omega}}\left( \frac{\pi}{4}-\frac{\pi }{\sqrt{3}}\frac{w_\Omega }{P(\Omega)}+ \frac{\pi }{3}\frac{w_\Omega^2 }{P(\Omega)^2}\right)\\
    &\geq \frac{w_\Omega^2}{\abs{\Omega}}\left( \frac{\pi}{4}-\frac{\pi }{2\sqrt{3}}\frac{w_\Omega }{\diam(\Omega)}\right)\\
     &\geq \frac{w_\Omega^2}{\abs{\Omega}}\left( \frac{\pi}{4}-\frac{\pi }{2\sqrt{3}}\frac{\sigma }{K(2)}\right)\\
     &\geq \frac{4}{3\sqrt{3}} \frac{w_\Omega^2}{\abs{\Omega}}
\end{split}
\end{equation}
Finally, by combining  \eqref{pr_vs_pw} and \eqref{pw-m}, we get the conclusion
\begin{equation}\label{final_final}
 \frac{T(\Omega)P^2(\Omega)}{|\Omega|^3}-\frac{1}{3} \ge \frac{1}{6}\left[\frac{P(\Omega)R_\Omega}{\abs{\Omega}} -1\right]^3\ge 
 \tilde{K}\left[\frac{\abs{Q\bigtriangleup \Omega}}{\abs{\Omega}}\right]^3.
\end{equation}

\end{proof}

The next remark shows that a sequence of thinning triangles is not sharp for \eqref{quantitative_width} in the case $n=2$ and this is the reason for which we need Theorem \ref{teorema3} to obtain more precise information. 

\begin{oss} \label{triangle} 

Let us consider a sequence of isosceles triangles $\mathcal{T}_l$ of base $L$ and height $l$ such that $ \abs{\mathcal{T}_l}=1$.
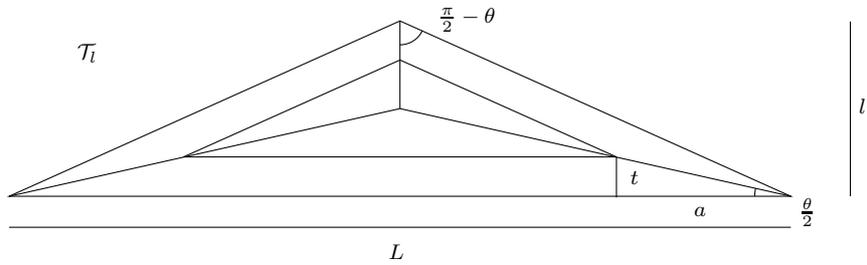
\begin{figure}[h]
\begin{center}
\begin{tikzpicture}[x=0.75pt,y=0.75pt,yscale=-0.6,xscale=0.6]

\draw   (325.5,203) -- (650.5,350) -- (0.5,350) -- cycle ;
\draw    (0.5,350) -- (325.5,276.5) ;
\draw    (325.5,276.5) -- (650.5,350) ;
\draw    (325.5,203) -- (325.5,276.5) ;
\draw    (145.5,317) -- (505.5,317) ;
\draw    (325.5,235.75) -- (505.5,317) ;
\draw    (325.5,235.75) -- (145.5,317) ;
\draw    (505.5,317) -- (505.5,351) ;
\draw    (700,203.5) -- (700,350) ;
\draw    (0.5,376) -- (650.5,376) ;
\draw  [draw opacity=0] (620.5,349.9) .. controls (620.51,347.64) and (620.76,345.44) .. (621.24,343.33) -- (650.5,350) -- cycle ; \draw   (620.5,349.9) .. controls (620.51,347.64) and (620.76,345.44) .. (621.24,343.33) ;
\draw  [draw opacity=0] (344.12,211.04) .. controls (340.97,218.08) and (333.82,223) .. (325.5,223) .. controls (325.41,223) and (325.31,223) .. (325.22,223) -- (325.5,203) -- cycle ;
\draw   (344.12,211.04) .. controls (340.97,218.08) and (333.82,223) .. (325.5,223) .. controls (325.41,223) and (325.31,223) .. (325.22,223) ;

\draw (515,327) node [anchor=north west][inner sep=0.75pt]   [align=left] {\footnotesize $t$};
\draw (705,266) node [anchor=north west][inner sep=0.75pt]   [align=left] {\footnotesize $l$};
\draw (314,388) node [anchor=north west][inner sep=0.75pt]   [align=left] {\footnotesize $L$};
\draw (568,357) node [anchor=north west][inner sep=0.75pt]   [align=left] {\footnotesize $a$};
\draw (655.22,349.67) node [anchor=north west][inner sep=0.75pt]   [align=left] {\footnotesize $\frac{\theta}{2}$};
\draw (355,189.67) node [anchor=north west][inner sep=0.75pt]   [align=left] {\footnotesize $\frac{\pi}{2}- \theta$};
\draw (55,220.67) node [anchor=north west][inner sep=0.75pt]   [align=left] {\footnotesize $\mathcal{T}_l$};
\end{tikzpicture}

\end{center}
\caption{Isosceles triangle $\mathcal{T}_l$ of base $L$ and height $l$.} \label{fig:M4}
\end{figure}

If we compute \eqref{quantitative_R} on the sequence $\mathcal{T}_l$ and we use \eqref{convex_estimates}, we get, for every $l$, 
\begin{align}
\label{fru}
    \dfrac{T(\mathcal{T}_l) P^2(\mathcal{T}_l)}{|\mathcal{T}_l|^3}-\frac{1}{3}\geq 
    \frac{1}{6}\left(\frac{P(\mathcal{T}_l)R_{\mathcal{T}_l}}{\abs{\mathcal{T}_l}}-1\right)^3=\frac{1}{6}
     \end{align}
     and, so,  the quantity on the left-hand side of  \eqref{fru} is bounded away from zero. 
\end{oss}

\begin{oss}
We point out that 
\begin{equation*}
    \frac{P(\Omega)R_\Omega}{\abs{\Omega}} -1\ge 
K \frac{\abs{Q\bigtriangleup \Omega}}{\abs{\Omega}},
\end{equation*}
in \eqref{final_final}  is a quantitative version of the inequality in the right hand side of \eqref{convex_estimates}.
\end{oss}

\begin{open} We conclude by listing the following open problems:
\begin{itemize}
    \item We believe that the exponent 3 in the inequality \eqref{quantitativa2}  
     is not sharp: we expect it to be 1. 
    We clarify that in  Example \eqref{esret}.

    \item 
    We conjecture that the sharp exponent in \eqref{quantitative_width} in the case $n>2$ is $1$ (see Remark \ref{remarkn}). 
    \item The results contained in Theorem \ref{teorema3} could be studied in higher dimension and extended 
       to the $f,p$-torsional rigidity.
    Our proof cannot be adapted to higher dimension because in dimension $n>2$ we do not have any more Steiner formulas for inner parallel sets \eqref{steiner1} and \eqref{steiner2}. 
\end{itemize}
\end{open}
\begin{esempio}\label{esret}
Let $\displaystyle{\Omega_l= \left(-\frac{1}{2l},\frac{1}{2l} \right)\times \left(-\frac{l}{2},\frac{l}{2}\right)}$ be a sequence of rectangles of measure $1$. It is possible to give an explicit upper bound to the functional $\mathcal{F}_2(\Omega_l)$. Hence, following the computations in \eqref{p-tortionup}, we have
$$
\mathcal{F}_2 (\Omega_l)- c_2 \leq 2 l^2.$$
Considering the rectangle $Q$ with sides $P(\Omega_l)/ 2$ and $w_{\Omega}$ containing $\Omega_l$, that is 
$$\displaystyle{Q= \left(-\frac{1+l^2}{2l},\frac{1+l^2}{2l} \right)\times \left(-\frac{l}{2},\frac{l}{2}\right)},$$
it is straightfoward to compute
$$\abs{\Omega_l \Delta Q}= 2l^2.$$
Hence, we have 
$$
2 \geq \frac{\mathcal{F}_2(\Omega_l) - c_2}{\abs{\Omega_l \Delta Q}} \geq \tilde{K}\,\abs{\Omega_l \Delta Q}^2 =\tilde{K}\, l^4.
$$
\end{esempio}

\section*{Acknowledgements}
This work has been partially supported  by GNAMPA of INdAM. In particular, the author Gloria Paoli is supported by the Alexander von Humboldt Foundation with an Alexander von Humboldt research fellowship.  Moreover, the authors would like to thank the reviewers for their suggestions to improve
this paper.

\bibliographystyle{plain}
\bibliography{biblio}

\Addresses 

\end{document}